\documentclass[reqno,12pt]{article} \usepackage{fontenc}
\usepackage{amsfonts}
\usepackage[utf8]{inputenc}
\usepackage[english]{babel}
\usepackage{amsmath}
\usepackage{amssymb}

\newtheorem{theorem}{Theorem}

\newtheorem{lemma}[theorem]{Lemma}

\newtheorem{example}{Example}[theorem]
\newtheorem{definition}{Definition}[theorem]
\begin{document}
{\small {~~\\ \vspace{-3.5cm} \begin{flushright}\jobname \\ \today
\end{flushright} \vspace{1.5cm} }}

\begin{center}
{\bf \Large
COMPLEX FUZZY LIE ALGEBRAS}\\~~\\
{Shadi Shaqaqha}\\~\\
Yarmouk University, Irbid, Jordan\\
shadi.s@yu.edu.jo
\end{center}

%
%
%

\begin{abstract}
A complex fuzzy Lie algebra is a fuzzy Lie algebra whose membership function takes values in the unit
circle in the complex plane. In this paper, we define the complex fuzzy Lie subalgebras and complex fuzzy ideals of
Lie algebras. Then, we investigate some of characteristics of complex fuzzy Lie subalgebras. The relationship between
complex fuzzy Lie subalgebras and fuzzy Lie subalgebras is also investigated. Finally, we define the image and the
inverse image of complex fuzzy Lie subalgebra under Lie algebra homomorphism. The properties of complex fuzzy Lie
subalgebras and complex fuzzy ideals under homomorphisms of Lie algebras are studied. \end{abstract}


\vspace{.4cm} \noindent {\bf Keywords:} Lie algebras; Lie homomorphism; $\pi$-fuzzy set; $\pi$-fuzzy Lie subalgebra; homogeneous
complex fuzzy set; complex fuzzy Lie subalgebra; complex fuzzy Lie ideal.
\section{INTRODUCTION}
Lie algebras were introduced by Sophus Lie (1842-1899) while he was attempting to classify certain smooth subgroups of general linear groups that are now called Lie groups. By now, both Lie groups and Lie algebras have become essential
to many parts of mathematics and physics.\\
The notion of fuzzy sets was firstly introduced by Zadeh \cite{Zadeh}. The fuzzy set theory states that there are propositions with an infinite number of truth values, assuming two extreme values, $1$ (totally true), $0$ (totally false) and a continuum in between, that justify
the term ”fuzzy”.  Applications of this
theory can be found, for example, in artificial intelligence, computer science, control engineering, decision
theory, logic and management science.\\
After the introduction of fuzzy sets by Zadeh, there are a number of generalizations of this fundamental concept. The
notion of complex fuzzy sets is introduced by Ramot and Milo \cite{Ramot} is one among them. Here, the membership values are
complex numbers, drawn from the unit disc of the complex plane. For more details on complex fuzzy sets, we refer the readers to \cite{Ramot1, Ramot}. Possible
applications, which demonstrate the complex fuzzy sets  concept, include a complex fuzzy
representation of solar activity (via measurements of the sunspot number),
signal processing application, time series forecasting problems and compare
the two national economy (\cite{Suresh}).
\\
The fuzzy Lie subalgebras and fuzzy Lie ideals are considered in \cite{Kim} by Kim and
Lee, and in \cite{Yehia1,Yehia2} by Yehia. They established the analogues of most of the fundamental ground results
involving Lie algebras in the fuzzy setting. Recently, the complex fuzzy subgroups and subrings are studied in \cite{Alsarahead, Alsarahead2} as a generalization of fuzzy groups and fuzzy rings, respectively.\\
In this paper we describe complex fuzzy Lie algebras.
\section{PRELIMINARIES}
In this section we give some relevant definitions, notations, and results that will be used frequently throughout the paper. Romot and Milo have defined complex fuzzy sets as follows (\cite{Ramot}).
\begin{definition}\label{CFL1}{
Let $X$ be a nonempty set. A {\em complex fuzzy set} $A$ on $X$ is an object having the form $A=\left\{(x, \mu_A(x))~|~x\in X\right\}$, where $\mu_A$ denotes the degree of membership function that assigns each element $x\in X$ a complex number $\mu_A(x)$ lies within the unit circle in the complex plane.
}
\end{definition}
We shall assume that is $\mu_A(x)$ will be represented by $r_A(x)e^{iw_A(x)}$ where $i=\sqrt{-1}, r_A(x)\in[0, 1]$, and $w_A(x)\in[0, 2\pi]$. Note that by setting $w_A(x)=0$ in the definition above, we return back to the traditional fuzzy subset.\\
Let $\mu_1=r_1e^{iw_1}$, and $\mu_2=r_2e^{iw_2}$ be two complex numbers lie within the unit circle in the complex plane. By $\mu_1\leq \mu_2$, we mean $r_1\leq r_2$ and $w_1\leq w_2$. A complex fuzzy set $A=\left\{(x, \mu(x))~|~x\in X \right\}$, with membership function $\mu(x)=r(x)e^{iw(x)}$, is homogeneous if for all $x, y\in X$ we have $r(x)\leq r(y)$ if and only if $w(x)\leq w(y)$. Let $A= \left\{x, \mu_A (x))~|~x\in X\right\}$ and $B= \left\{x, \mu_B (x))~|~x\in X\right\}$ be two complex fuzzy subsets of the same set $X$, with membership functions $\mu_A(x)=r_A(x)e^{iw_A(x)}$ and $\mu_B(x)=r_B(x)e^{iw_B(x)}$, respectively. Then, $A$ is said to be homogeneous with $B$ if for all $x, y\in X$, we have $r_A(x)\leq r_B(y)$ if and only if $w_A(x)\leq w_B(y)$. Let $A= \left\{x, \mu (x))~|~x\in X\right\}$ be a homogenous complex fuzzy set. One can easily show if $x, y\in X$, then either we have $\mu_A (x)\leq \mu_A (y)$ or $\mu_A (y)\leq \mu_A (x)$.\\
In \cite{Alsarahead}, the authors gave the following definition.
\begin{definition}\label{CFL1A}
Let $A=\left\{(x, \mu_A(x)~|~x\in X\right\}$ be  fuzzy set. Then the set $A_{\pi}=\left\{(x, \gamma_A(x)~|~x\in X\right\}$ is said to be a {\em $\pi$-fuzzy subset} where $\gamma_A(x)=2\pi\mu_A(x)$.
\end{definition}
Let $A= \left\{x, \mu_A (x))~|~x\in X\right\}$ be a complex fuzzy set. For the sake of simplicity, we shall use the notation $A= \mu_A$. Also, throughout this paper, all complex fuzzy sets are homogeneous.\\
Let $F$ be a ground field. A Lie algebra over $F$ is a vector space $L$, together with a $F$-bilinear mapping
$$L\times L\rightarrow L : (x,y)\mapsto [x,y]$$
called a Lie product (or Lie bracket), and satisfying the following two identities, for any $x, y, z \in L$:\\
\begin{itemize}
\item[(i)] $[x,x]=0$ (this identity implies antisymmetry, that is $[x,y]=-[y,x]$).
\item[(ii)] $[x, [y,z]]+ [y,[z,x]]+[z,[x,y]]=0$ (Jacobi identity).
\end{itemize}
Subalgebras of Lie algebras and homomorphisms (or isomorphisms) between Lie algebras are defined as usual (\cite{Bahturin 1}).\\
Throughout this paper, $L$ is a Lie algebra over field $F$.
\begin{definition}\label{CFL2A}
Let $A_{\pi}=\left\{(x, \gamma_A(x)~|~x\in L\right\}$ be a $\pi$-fuzzy subset of $L$. Then $A_{\pi}$ is said to be {\em $\pi$-fuzzy subalgebra} of $L$ if the following conditions are satisfied for all $x, y\in L$ and $\alpha\in F$:
\begin{itemize}
\item[(i)] $\gamma_A(x+y)\geq \gamma_A(x)\wedge \gamma_A(y)$,
\item[(ii)] $\gamma_A(\alpha x)\geq \gamma_A(x)$,
\item[(iii)] $\gamma_A([x, y])\geq\gamma_A(x)\wedge \gamma_A(y)$.
\end{itemize}
\end{definition}
If the condition (iii) is replaced by $\gamma_A([x, y])\geq \gamma_A(x) \vee\gamma_A(y)$, then $\gamma_A$ is called a {\em $\pi$-fuzzy ideal} of $L$.\\
The following facts are obvious. The proofs are very similar to the proof of analogous results for fuzzy subgroups \cite{Alsarahead} and fuzzy subrings \cite{Alsarahead2}.
\begin{theorem}\label{CFL3A}
A $\pi$-fuzzy set $A_{\pi}$ is a $\pi$-fuzzy subalgebra of $L$ if and only if $A$ is a fuzzy subalgebra of $L$.
\end{theorem}
\begin{theorem}\label{CFL4A}
A $\pi$-fuzzy set $A_{\pi}$ is a $\pi$-fuzzy ideal of $L$ if and only if $A$ is a fuzzy ideal of $L$.
\end{theorem}
\section{Complex Fuzzy Lie Algebras}
For the sake of simplicity, we shall use the symbols $a\wedge b =\mathrm{min}\left\{a, b\right\}$ and $a\vee b =\mathrm{max}\left\{a, b\right\}$.
\begin{definition}\label{CFL2}
A (homogeneous) complex fuzzy set $\mu_A$ on $L$ is a {\em complex fuzzy subalgebra} if the following conditions are satisfied for all $x, y\in L$, and $\alpha\in F$:
\begin{itemize}
\item[(i)] $\mu_A(x+y)\geq \mu_A(x)\wedge \mu_A(y)$,
\item[(ii)] $\mu_A(\alpha x)\geq \mu_A(x)$,
\item[(iii)] $\mu_A([x, y])\geq\mu_A(x)\wedge \mu_A(y)$.
\end{itemize}
\end{definition}
If the condition (iii) is replaced by $\mu_A([x, y])\geq \mu_A(x) \vee\mu_A(y)$, then $\mu_A$ is called a {\em complex fuzzy ideal} of $L$. Note that the second condition implies $\mu_A(x)\leq \mu_A(0)$ and $\mu_A(-x)\leq\mu_A(x)$ for all $x\in L$.\\
It is clear that if $\mu_A$ is a complex fuzzy ideal of $L$, then it is a complex fuzzy subalgebra of $L$. But the converse is not true in general as can be seen in the following example.
\begin{example}\label{CFL1B}
It is well known that the set $\mathbb{R}^3=\{(x, y, z)~|~x, y, z\in\mathbb{R}\}$ of all $3$-dimensional real vectors forms a Lie algebra over $F=\mathbb{R}$ with the usual cross product $\times$. Define $A=\mu_A$, where $\mu_A: \mathbb{R}^3\rightarrow E^2$ ($E^2$ is the unit disc), by
$$\mu_{A}(x)=\left\{
  \begin{array}{lr}
  0.9e^{i\frac{3\pi}{2}} &: x=y=z=0 \\
   0.6e^{i\frac{\pi}{2}} &: x\neq 0, y=z=0\\
	0 &: otherwise
 \end{array}
 \right.$$
Then it is clear that $A$ is a complex fuzzy subalgebra of $L=\mathbb{R}^3$. But it is not a complex fuzzy Lie ideal since $\mu_A([(1, 0, 0), (1, 1, 1)])=\mu_A(0, -1, 1)=0\ngeq \mu_A(1, 0, 0)$.
\end{example}
The following lemma is obvious. The proof is very similar to the proof of analogous result for the traditional fuzzy Lie algebra.
\begin{lemma}\label{CFL1AA}
Let $A=\mu_A$ be a complex fuzzy subalgebra of $L$. Then
\begin{itemize}
\item[(i)] $\mu_A(-x)=\mu_A(x)$ for all $x\in L$,
\item[(ii)] $\mu_A(x-y)= \mu_A(0)$ implies $\mu_A(x)=\mu_A(y)$,
\item[(iii)] $\mu_A(x)<\mu_A(y)$ implies $\mu_A(x-y)=\mu_A(x)=\mu_A(y-x)$.
\end{itemize}
\end{lemma}
The following theorem gives the relationship between complex fuzzy subalgebras and fuzzy subalgebras. There is an analogous result for groups \cite{Alsarahead} and rings \cite{Alsarahead2}.
\begin{theorem}\label{CLF11A}
Let $A=\{(x, \mu_A(x))~|~x\in L\}$ be a complex fuzzy set with membership function $\mu_A(x)=r_A(xe^{iw_A(x)}$. Then $A$ is a complex fuzzy subalgebra if and only if the fuzzy subset $\overline{A}=\{(x, r_A(x))~|~x\in L\}$ is a fuzzy subalgebra and the $\pi$-fuzzy set $\underline{A}=\{(x, w_A(x))~|~x\in L\}$ is a $\pi$-fuzzy subalgebra.
\end{theorem}
{\it Proof.~}
Let $A$ be a complex fuzzy subalgebra of $L$. For $x, y\in L$ and $\alpha \in F$, we have
\begin{eqnarray*}
r_A(x+y)e^{i(w_A(x+y)}&=& \mu_A(x+y)\\
&\geq & \mu_A(x)\wedge \mu_A(y)\\
&=& \left(r_A(x)\wedge r_A(y)\right)e^{i(w_A(x)\wedge w_A(y))}~~~(A~is~\mathrm{homogeneous}).
\end{eqnarray*}
Thus, $r_A(x+y)\geq r_A(x)\wedge r_A(y)$ and $w_A(x+y)\geq w_A(x)\wedge w_A(y)$. The proof of $r_A(\alpha x)\geq r_A(x)$, $w_A(\alpha x)\geq w_A(x)$, $r_A([x, y])\geq r_A(x)\wedge r_A(y)$, and $w_A([x, y])\geq w_A(x)\wedge w_A(y)$. are very similar to the previous proof, so we omit them. Hence $\overline{A}$ is a fuzzy subalgebra and $\underline{A}$ is a $\pi$-fuzzy subalgebra.\\
Conversely, suppose that $\overline{A}$ and $\underline{A}$ are fuzzy subalgebra and $\pi$-fuzzy subalgebra, respectively, of $L$. Then for any $x, y\in L$ and $\alpha\in F$, we have
\begin{eqnarray*}
r_A(x+y)&\geq & r_A(x)\wedge r_A(y),\\
r_A(\alpha x)&\geq& r_A(x),\\
\mathrm{and}~r_A([x, y])&\geq& r_A(x)\wedge r_A(y).
\end{eqnarray*}
Also,
\begin{eqnarray*}
w_A(x+y)&\geq & w_A(x)\wedge w_A(y),\\
w_A(\alpha x)&\geq& w_A(x),\\
\mathrm{and}~w_A([x, y])&\geq& w_A(x)\wedge w_A(y).
\end{eqnarray*}
Now for any $x, y\in L$ and $\alpha\in F$, we have
\begin{eqnarray*}
\mu_A(x+y)&=& r_A(x+y)e^{i(w_A(x+y))}\\
&\geq& (r_A(x)\wedge r_A(y))e^{i(w)A(x)\wedge w_A(y))}\\
&=&r_A(x)e^{iw_A(x)}\wedge r_A(y)e^{iw_A(y)}\\
&=&\mu_A(x)\wedge\mu_A(y).
\end{eqnarray*}
Similarly, one can prove $\mu_A(\alpha x)\geq \mu_A(x)$ and $\mu_A([x, y])\leq \mu_A(x)\wedge\mu_A(y)$.
This shows $A$ is a complex fuzzy subalgebra of $L$.\hfill $\Box$

In a similar way, one can show the following theorem.
\begin{theorem}\label{CLF12A}
Let $A=\{(x, \mu_A(x))~|~x\in L\}$ be a complex fuzzy set with membership function $\mu_A(x)=r_A(xe^{iw_A(x)}$. Then $A$ is a complex fuzzy ideal if and only if the fuzzy subset $\overline{A}=\{(x, r_A(x))~|~x\in L\}$ is a fuzzy ideal and the $\pi$-fuzzy set $\underline{A}=\{(x, w_A(x))~|~x\in L\}$ is a $\pi$-fuzzy ideal.
\end{theorem}
Let $V$ be a vector space. For $t$ lies in the unit disk in the complex plane and complex fuzzy set $\mu_A$, the set $U(\mu_A, t)= \left\{x\in V~|~\mu_A(x)\geq t\right\}$ is called an upper level of $\mu_A$. In particular, if $t=\alpha e^{i0}; \alpha\in [0, 1]$, we get the upper level subset $A_{\alpha}=\{x\in X~|~r_A(x)\geq \alpha\}$. Also, if $t=0e^{i\beta}$ ($\beta\in [0, 2\pi]$), we get the level subset $A_{\beta}=\{x\in X~|~w_A(x)\geq \beta\}$. The following theorem will show the relations between complex fuzzy Lie subalgebras of $L$ and Lie subalgebras of $L$.
\begin{theorem}\label{CFL3}
Let $\mu_A$ be a complex fuzzy subset of $L$. Then the following statements are equivalent:
\begin{itemize}
\item[(i)] $\mu_A$ is a complex fuzzy subalgebra of $L$,
\item[(ii)] the upper level $U(\mu_A, t)$ is a subalgebra of $L$ for every $t\in \mathrm{Im}(\mu_A)$.
\end{itemize}
\end{theorem}
{\it Proof.~}
For $t\in\mathrm{Im}(\mu_A)$, let $x, y\in U(\mu_A, t)$, and $\alpha\in F$. As $\mu_A$ is a complex fuzzy subalgebra of $L$, we have $\mu_A(x+y)\geq\mu_A(x)\wedge\mu_A(y)\geq t$, $\mu_A(\alpha x)\geq\mu_A(x)\geq t$, and $\mu_A([x, y])\geq\mu_A(x)\wedge\mu_A(y)\geq t$, and so $x+y$, $\alpha x$, and $[x, y]$ are elements in $U(\mu_A, t)$. Conversely, suppose that the upper levels $U(\mu_A, t)$ are Lie subalgebras of $L$ for every $t\in \mathrm{Im}(\mu_A)$. Let $x, y\in L$ and $\alpha\in F$. We may assume $\mu_A(y)\geq \mu_A(x)=t_1$, so $x, y\in U(\mu_A, t_1)$. Since $U(\mu_A, t_1)$ is a subspace of $L$, we have $x+y$ and $\alpha x$ are in $U(\mu_A, t_1)$, and so $\mu_A(\alpha x)\geq t_1= \mu_A(x)$ and $\mu_A(x+y)\geq t_1 =\mu_A(x)\wedge\mu_A(y)$. Also, since $U(\mu_A, t_1)$ is a Lie subalgebra of $L$, we have $[x, y]\in U(\mu_A, t_1)$. Hence, $\mu_A([x, y])\geq t_1 =\mu_A(x)\wedge\mu_A(y)$.\hfill $\Box$

Similarly, one can prove the following theorem.
\begin{theorem}\label{CFL4}
Let $\mu_A$ be a complex fuzzy subset of $L$. Then the following statements are equivalent:
\begin{itemize}
\item[(i)] $\mu_A$ is a complex fuzzy ideal of $L$,
\item[(ii)] the upper level $U(\mu_A, t)$ is an ideal of $L$ for every $t\in \mathrm{Im}(\mu_A)$.
\end{itemize}
\end{theorem}
Let $V$ be a vector space. For $t$ lies in the unit disk in the complex plane and complex fuzzy set $\mu_A$, the set $U(\mu_A^>, t)= \left\{x\in V~|~\mu_A(x)> t\right\}$ is called a {\em strong upper level} of $\mu_A$. We have the following result (there is a similar result in the case of polygroups \cite{NON}).
\begin{theorem}\label{CFL3B}
Let $\mu_A$ be a complex fuzzy subset of $L$. Then the following statements are equivalent:
\begin{itemize}
\item[(i)] $\mu_A$ is a complex fuzzy subalgebra of $L$,
\item[(ii)] the strong upper level $U(\mu_A^>, t)$ is a subalgebra of $L$ for every $t\in \mathrm{Im}(\mu_A)$.
\end{itemize}
\end{theorem}
{\it Proof.~}
For $t\in\mathrm{Im}(\mu_A)$, let $x, y\in U(\mu_A^>, t)$, and $\alpha\in F$. As $\mu_A$ is a complex fuzzy subalgebra of $L$, we have $\mu_A(x+y)\geq\mu_A(x)\wedge\mu_A(y)>t$, $\mu_A(\alpha x)>\mu_A(x)\geq t$, and $\mu_A([x, y])\leq\mu_A(x)\wedge\mu_A(y)> t$, and so $x+y$, $\alpha x$, and $[x, y]$ are elements in $U(\mu_A^>, t)$. Conversely, assume that every strong upper level $U(\mu_A^>, t)$ is a Lie subalgebra of $L$ for every $t\in \mathrm{Im}(\mu_A)$. Let $x, y\in L$ and $\alpha\in F$. We need to show that the conditions of Definition \ref{CFL2} are satisfied. If $\mu_A(x)=0$ or $\mu_A(y)=0$, then the proof is obvious, so we may assume that $\mu_A(x)\neq 0$ and $\mu_A(y)\neq 0$. Let $t_0$ be the largest complex number on the closed unit disc of the complex plane such that $t_0<\mu_A(x)\wedge\mu_A(y)$ and there is no $a\in L$ satisfying $t_0< \mu_A(a)< \mu_A(x)\wedge\mu_A(y)$. Having $x, y\in U(\mu_A^>, t_0)$ implies that $x+y, [x, y]\in U(\mu_A^>, t_0)$, and hence $\mu_A(x+y), \mu_A([x, y])> t_0$. Since there exist no $a\in L$ with $t_0< \mu_A(a)< \mu_A(x)\wedge \mu_A(y)$, it follows that $\mu_A(x+y), \mu_A([x, y])>\mu_A(x)\wedge \mu_A(y)$. Again let $t_0$ be the largest complex number on the closed unit disc of the complex plane such that $t_0<\mu_A(x)$ and there is no $a\in L$ with $t_0<\mu_A(a)<\mu_A(x)$. As $U(\mu_A^>, t_0)$ is subalgebra, we have $\alpha x\in U(\mu_A^>, t_0)$, and so $\mu_A(\alpha x)> t_0$. Thus $\mu_A(\alpha x)> \mu_A(x)$. \hfill $\Box$

Similarly one can prove the following theorem.
\begin{theorem}\label{CFL3C}
Let $\mu_A$ be a complex fuzzy subset of $L$. Then the following statements are equivalent:
\begin{itemize}
\item[(i)] $\mu_A$ is a complex fuzzy ideal of $L$,
\item[(ii)] every strong upper level $U(\mu_A^>, t)$ is an ideal of $L$ for every $t\in \mathrm{Im}(\mu_A)$.
\end{itemize}
\end{theorem}

Let $\mu_A$ and $\mu_B$ be complex fuzzy subsets of $L$. Let us introduce the complex fuzzy sum of $\mu_A$ and $\mu_B$ of $L$ as follows:
$$\mu_{A+B}(x)= \mathrm{sup}_{x=a+b}\left\{\mu_A(a)\wedge\mu_B(b)\right\}.$$
Now, we have the following theorem.
\begin{theorem}\label{CFL5}{
Let $A$ and $B$ be (homogenous) complex fuzzy ideals of $L$ with membership functions $\mu_A(x)=r_A(x)e^{iw_A(x)}$ and $\mu_B(x)=r_B(x)e^{iw_B(x)}$, respectively such that $A$ is homogenous with $B$. Then so $\mu_{A+B}$ is a complex fuzzy ideal too.
}
\end{theorem}
{\it Proof.~}
First we show that $A+B$ is homogeneous. For $x\in L$, we have
\begin{eqnarray*}
\mu_{A+B}(x)&=&\mathrm{sup}_{x=a+b}\{\mu_A(a)\wedge\mu_B(b)\}\\
&=& \mathrm{sup}_{x=a+b}\{r_A(a)e^{iw_A(a)}\wedge r_B(b)e^{iw_B(b)}\}\\
&=&\mathrm{sup}_{x=a+b}\{r_A(a)\wedge r_B(b)e^{i(w_A(a)\wedge w_B(b))}\}\\
&=&\mathrm{sup}_{x=a+b}\{r_A(a)\wedge r_B(b)\}e^{i\mathrm{sup}_{x=a+b}\{w_A(a)\wedge w_B(b)\}}.
\end{eqnarray*}
Let $x, y\in L$ such that $\mathrm{sup}_{x=a+b}\{r_A(a)\wedge r_B(b)\}\leq \mathrm{sup}_{y=c+d}\{r_A(c)\wedge r_B(d)\}$. Suppose that $\mathrm{sup}_{y=c+d}\{w_A(c)\wedge w_B(d)\}<\mathrm{sup}_{x=a+b}\{w_A(a)\wedge w_B(b)\}$. Then, we can find $a_1, b_1\in L$ such that $x=a_1+ b_1$ and $\mathrm{sup}_{y=c+d}\{w_A(c)\wedge w_B(d)\}< w_A(a_1)\wedge w_B(b_1)$. Without loss of generality, we may assume $w_A(a_1)\leq \mu_B(b_1)$. Hence if $c, d\in L$ with $y=c+d$, then $w_A(c)\wedge w_B(d)< w_A(a_1)$. If $w_A(c)< w_B(d)$, then $w_A(c)< w_A(a_1)$, and from the homogeneity of $A$, we have $r_A(c)< r_A(a_1)$. Also if $w_B(d)< w_A(c)$, then $w_B(d)< w_A(a_1)$. Since $A$ is homogenous with $B$, it follows $r_B(d)< r_A(a_1)$. Thus $\mathrm{sup}_{y=c+d}\{r_A(c)\wedge r_B(b)\}< r_A(a_1)$, and so $\mathrm{sup}_{x=a+b}\{r_A(a)\wedge r_B(b)\}>\mathrm{sup}_{y=c+d}\{r_A(c)\wedge r_B(d)\}$. Contradiction.\\
Let $x, y\in L$. Then
\begin{eqnarray*}
\mu_{A+B}(x)\wedge\mu_{A+B}(y) &=& \mathrm{sup}_{x=a+b}\left\{\mu_A(a)\wedge \mu_B(b)\right\}\wedge \mathrm{sup}_{y=c+d}\left\{\mu_A(c)\wedge\mu_B(d)\right\}\\
&=&\mathrm{sup}_{x=a+b, y=c+d}\left\{(\mu_A(a)\wedge \mu_B(b))\wedge (\mu_A(c)\wedge\mu_B(d))\right\}\\
&=&\mathrm{sup}_{x=a+b, y=c+d}\left\{(\mu_A(a)\wedge \mu_A(c))\wedge (\mu_B(b)\wedge\mu_B(d))\right\}\\
&\leq&\mathrm{sup}_{x=a+b, y=c+d}\left\{\mu_A(a+c)\wedge \mu_B(b+d)\right\}\\
&=&\mu_{A+B}(x+y).
\end{eqnarray*}
Also, for $\alpha\in F$ and $x\in L$ with $\mu_A(x)\neq 0$, we have
\begin{eqnarray*}
\mu_{A+B}(x)&=&\mathrm{sup}_{x=a+b}\left\{\mu_A(a)\wedge \mu_B(b)\right\}\\
&\leq& \mathrm{sup}_{x=a+b}\left\{\mu_A(\alpha a)\wedge \mu_B(\alpha b)\right\}\\
&=&\mathrm{sup}_{\alpha x=\alpha a+ \alpha b}\left\{\mu_A(\alpha a)\wedge \mu_B(\alpha b)\right\}\\
&\leq&\mathrm{sup}_{\alpha x =c+d}\left\{\mu_A(c)\wedge \mu_B(d)\right\}\\
&=&\mu_{A+B}(\alpha x).
\end{eqnarray*}
Let $x, y\in L$, we need to show that $\mu_{A+B}([x, y])\geq \mu_{A+b}(x)\vee \mu_{A+B}(y)$. Suppose that $\mu_{A+B}([x, y])< \mu_{A+B}(x)\vee \mu_{A+B}(y)$. Then $\mu_{A+B}([x, y])< \mu_{A+B}(x)$ or $\mu_{A+B}([x, y])< \mu_{A+B}(y)$. Without loss of generality, we may assume $\mu_{A+B}([x, y)< \mu_{A+B}(x)$. Choose $t$ within the unit circle in the complex plane with $\mu_{A+B}([x, y])< t< \mu_{A+B}(x)$.  Thus there exist $a_1, b_1\in L$ with $x=a_1+b_1$, $\mu_A(a_1)> t$, and $\mu_B(b_1)>t$. Hence,
\begin{eqnarray*}
\mu_{A+B}([x, y]&=&\mu_{A+B}([a_1+b_1, y])\\
&>& \mathrm{sup}_{[x, y]= [a, y]+[b,y]}\left\{\mu_{A}([a, y], \mu_B([b, y])\right\}\\
&>&\left(\mu_A([a_1, y])\wedge \mu_B([b_1, y])\right)\\
&\geq& (\mu_A(a_1)\vee \mu_A(y))\wedge (\mu_B(b_1)\vee \mu_B(y))\\
&\geq& t\\
&>&\mu_{A+B}([x, y]).
\end{eqnarray*}
 Contradiction. This shows $\mu_{A+B}$ is a complex fuzzy ideal of $L$
\hfill $\Box$

Recall that if $A= \left\{x, \mu_A (x))~|~x\in X\right\}$ and $B= \left\{x, \mu_B (x))~|~x\in X\right\}$ are two complex fuzzy subsets of the same set $X$, with membership functions $\mu_A(x)=r_A(x)e^{iw_A(x)}$ and $\mu_B(x)=r_B(x)e^{iw_B(x)}$, respectively, then $A\cap B=\{(x, \mu_{A\cap B}(x))~|~x\in X\}$, where
$$\mu_{A\cap B}(x)=r_{A\cap B}(x)e^{iw_{A\cap B}(x)}=(r_A(x)\wedge r_B(x))e^{i(w_A(x)\wedge w_B(x))}.$$
We have the following result.
\begin{theorem}\label{CFL6}
Let $A$ and $B$ be (homogenous) complex fuzzy subsets of $L$ such that $A$ is homogenous with $B$. If $A$ and $B$ are complex fuzzy subalgebras of $L$, then $A\cap B$ is a homogenous complex fuzzy subalgebra of $L$.
\end{theorem}
{\it Proof.~}
First, we need to show that $A\cap B$ is a homogenous complex fuzzy subset of $L$. Let $x, y\in L$ with $r_{A\cap B}(x)\leq r_{A\cap B}(y)$. Without loss of generality, we may assume $r_A(x)\leq r_B(x)$. Then $r_{A\cap B}(x)=r_A(x)$, and so $r_A(x)\leq r_A(y)$ and $r_A(x)\leq r_B(y)$. Also, as $A$ is homogenous, we have $w_A(x)\leq w_A(y)$. From the homogeneity of $A$ in $B$, we have $w_A(x)\leq w_B(x)$ and $w_A(x)\leq w_B(y)$. Hence, $w_{A\cap B}(x)=w_A(x)= w_A(x)\wedge w_B(x)\leq w_{A\cap B}(y)=w_A(y)\wedge w_B(y)$. In a similar way, one can prove if $w_{A\cap B}(x)\leq w_{A\cap B}(y)$ we have $r_{A\cap B}(x)\leq r_{A\cap B}(y)$.\\

Let $x, y\in L$ and $\alpha\in F$. We need to show that the conditions of Definition \ref{CFL2} are satisfied. From the homogeneity of $A\cap B$, it is enough to show
\begin{itemize}
\item[(i)] $r_{A\cap B}(x+y)\geq r_{A\cap B}(x)\wedge r_{A\cap B}(y)$,
\item[(ii)] $r_{A\cap B}(\alpha x)\geq r_{A\cap B}(x)$,
\item[(iii)] and $r_{A\cap B}([x, y]\geq r_{A\cap B}(x)\wedge r_{A\cap B}(y)$.
\end{itemize}
Now,
\begin{eqnarray*}
r_{A\cap B}(x+y)&=& r_A(x+y)\wedge r_B(x+y)\\
&\geq & (r_A(x)\wedge r_A(y))\wedge (r_B(x)\wedge r_B(y))~~~(\mathrm{Theorem~\ref{CLF11A}})\\
&=&(r_A(x)\wedge r_B(x))\wedge (r_A(y)\wedge r_B(y))\\
&=&r_{A\cap B}(x)\wedge r_{A\cap B}(y).
\end{eqnarray*}
Also
\begin{eqnarray*}
r_{A\cap B}(\alpha x)&=& r_A(\alpha x)\wedge r_B(\alpha x)\\
&\geq & r_A(x)\wedge r_B(y)~~~(\mathrm{Theorem~\ref{CLF11A}})\\
&=& r_{A\cap B}(x).
\end{eqnarray*}
In addition,
\begin{eqnarray*}
r_{A\cap B}([x, y])&=& r_A([x, y])\wedge r_B([x, y])\\
&\geq & (r_A(x)\wedge r_A(y))\wedge (r_B(x)\wedge r_B(y))~~~(\mathrm{Theorem~\ref{CLF11A}})\\
&=&(r_A(x)\wedge r_B(x))\wedge (r_A(y)\wedge r_B(y))\\
&=&r_{A\cap B}(x)\wedge r_{A\cap B}(y).
\end{eqnarray*}
\hfill $\Box$

Now, we can easily prove the following theorem too.
\begin{theorem}\label{CFL7}
Let $\{A_i~|~i\in I\}$ be a collection of complex fuzzy subalgebras of $L$ such that $A_j$ is homogenous with $A_k$ for all $j, k\in I$. Then $\bigcap_{i\in I}A_i=\mu_{\bigcap_{i\in I}A_i}$ is a complex fuzzy subalgebra of $L$, where
$$\mu_{\bigcap_{i\in I}A_i}=\left(\bigwedge_{i\in I}r_{A_i}\right)e^{i\left(\bigwedge_{i\in I}w_{A_i}\right)}.$$
\end{theorem}
We omit the proof of the following result since it is similar to the proof of above theorem.
\begin{theorem}\label{CFL8}
Let $\{A_i~|~i\in I\}$ be a collection of complex fuzzy ideals of $L$ such that $A_j$ is homogenous with $A_k$ for all $j, k\in I$. Then $\bigcap_{i\in I}A_i=\mu_{\bigcap_{i\in I}A_i}$ is a complex fuzzy ideal of $L$.
\end{theorem}
\section{On Lie algebra homomorphisms}
Let $f: X\rightarrow Y$ be a function. If $B=\mu_B$ is a complex fuzzy set of $Y$, then the preimage of $B$ is defined to be a complex fuzzy set $f^{-1}(B)=\mu_{f^{-1}}(B)$, where $\mu_{f^{-1}(B)}(x)=\mu_B(f(x))$ for any $x\in X$. Also if $A=\mu_A$ is a complex fuzzy set on $X$, then the image of $A$ is defined to be a complex fuzzy set $f(A)=\mu_{f(A)}$ where
$$\mu_{f(A)}(y)=\left\{
    \begin{array}{lr}
    \mathrm{sup}_{x\in f^{-1}(y)}\left\{\mu_A(x)\right\} &: y\in f(A) \\
     0 &: y\notin f(A).
 \end{array}
  \right.$$
(see for example \cite{Alsarahead}). The following theorem was obtained by  Alsarahead and Ahmad \cite{Alsarahead, Alsarahead2} in the setting of groups and rings.
\begin{theorem}\label{CFL9}
Let $\varphi: L\rightarrow L'$ be a Lie algebra homomorphism. If $B=\mu_B$ is a complex fuzzy Lie subalgebra of $L'$ with a membership function $\mu_{B}(y)=r_B(y)e^{iw_B(y)}$, then the complex fuzzy set $\varphi^{-1}(B)$ is also a complex fuzzy Lie subalgebra of $L$.
\end{theorem}
{\it Proof.~}
First, we need to show that $\varphi^{-1}(B)$ is homogenous. Note that if $x\in L$, then $\mu_{\varphi^{-1}(B)}(x)= \mu_B(\varphi(x))=r_B(\varphi(x))e^{iw_B(\varphi(x))}=(r_B\varphi(x))e^{iw_B(\varphi(x))}$. Now, if $x_1, x_2\in L$ with $(r_B\varphi)(x_1)\leq (r_B\varphi)(x_2)$; that is $r_B(\varphi(x_1))\leq r_B(\varphi(x_2))$, then from the homogeneity of $B$, we have $w_B(\varphi(x_1))\leq w_B(\varphi(x_2))$, and so $(w_B\varphi)(x_1)\leq (w_B\varphi)(x_2)$. Thus shows $\varphi^{-1}(B)$ is homogenous.\\
Let $x_1, x_2\in L$ and $\alpha\in F$. Then
\begin{eqnarray*}
\mu_{\varphi^{-1}(B)}(x_1+x_2)&=&\mu_B(\varphi(x_1+x_2))\\
&=& \mu_B(\varphi(x_1)+\varphi(x_2))~~~(\varphi~\mathrm{is~linear})\\
&\geq & \mu_B(\varphi(x_1))\wedge \mu_B(\varphi(x_2))\\
&=&\mu_{\varphi^{-1}(B)}(x_1)\wedge \mu_{\varphi^{-1}(B)}(x_2).
\end{eqnarray*}
Also,
\begin{eqnarray*}
\mu_{\varphi^{-1}(B)}(\alpha x_1)&=&\mu_B(\varphi(\alpha x_1))\\
&=& \mu_B(\alpha \varphi(x_1))~~~(\varphi~\mathrm{is~linear})\\
&\geq & \mu_B(\varphi(x_1))\\
&=&\mu_{\varphi^{-1}(B)}(x_1).
\end{eqnarray*}
Finally,
\begin{eqnarray*}
\mu_{\varphi^{-1}(B)}([x_1,x_2])&=&\mu_B(\varphi([x_1, x_2])\\
&=& \mu_B([\varphi(x_1), \varphi(x_2)])~~~(\varphi~\mathrm{is~homomorphism})\\
&\geq & \mu_B(\varphi(x_1))\wedge \mu_B(\varphi(x_2))\\
&=&\mu_{\varphi^{-1}(B)}(x_1)\wedge \mu_{\varphi^{-1}(B)}(x_2).
\end{eqnarray*}
\hfill $\Box$

\begin{theorem}\label{CFL9A}
Let $\varphi: L\rightarrow L'$ be a Lie algebra homomorphism. If $B=\mu_B$ is a complex fuzzy ideal of $L'$ with a membership function $\mu_{B}(y)=r_B(y)e^{iw_B(y)}$, then the complex fuzzy set $\varphi^{-1}(B)$ is also a complex fuzzy ideal of $L$.
\end{theorem}
{\it Proof.~}
This is very similar to the proof of Theorem \ref{CFL9}, and we shall omit it. \hfill $\Box$

It is known if $\varphi:L\rightarrow L'$ is a Lie algebra homomorphism and $A=\mu_A$ is a fuzzy subalgebra of $L$, then the image of $A$, $\varphi(A)$ is a fuzzy subalgebra of $L'$ (\cite{Kim}). We have the following result.
\begin{theorem}\label{CFL10}
Let $\varphi:L\rightarrow L'$ be a surjective Lie algebra homomorphism. If $A=\mu_A$, where $\mu_A(x)=r_A(x)e^{iw_A(x)}$ for any $x\in L$, is a complex fuzzy Lie subalgebra of $L$, then $\varphi(A)$ is also a complex fuzzy Lie subalgebra of $L'$.
\end{theorem}
{\it Proof.~}
We prove that $\varphi(A)$ is homogenous. Suppose $y\in L'$. Then
\begin{eqnarray*}
\mu_{\varphi(A)}(y)&=&\mathrm{sup}_{x\in \varphi^{-1}\{y\}}\{\mu_A(x)\}\\
&=&\mathrm{sup}_{x\in \varphi^{-1}\{y\}}\{r_A(x)e^{iw_A(x)}\}\\
&=& \left(\mathrm{sup}_{x\in \varphi^{-1}\{y\}}\{r_A(x)\}\right)e^{i\left(\mathrm{sup}_{x\in \varphi^{-1}\{y\}}\{w_A(x)\}\right)}\\
&=&{r_A}_{\varphi(A)}(y)e^{i{w_A}_{\varphi(A)}(y)}.
\end{eqnarray*}
Now let $y_1, y_2\in L'$ with ${r_A}_{\varphi(A)}(y_1)\leq{r_A}_{\varphi(A)}(y_2)$. Suppose that ${w_A}_{\varphi(A)}(y_2)<{w_A}_{\varphi(A)}(y_1)$. Then there exists $x_1\in \varphi^{-1}(\{y_1\})$ with ${w_A}_{\varphi(A)}(y_2)<{w_A}(x_1)$. Therefore, If $x\in \varphi^{-1}(\{y_2\})$, then $w_A(x)<{w_A}(x_1)$, and so, from the homogeneity of $A$, we obtain $r_A(x)< r_A(x_1)$. Thus, $\mathrm{sup}_{x\in \varphi^{-1}(\{y_2\})}\{r_A(x)\}< r_A(x_1)$, and so  ${r_A}_{\varphi(A)}(y_2)<{r_A}_{\varphi(A)}(y_1)$. Contradiction. This shows $\varphi(A)$ is homogenous.\\
 Since $A$ is a complex fuzzy subalgebra, it follows from Theorem \ref{CLF11A} that $\overline{A}=\{(x, r_A(x))~|~x\in L\}$ is a fuzzy subalgebra of $L$, and so the image ${r_A}_{\varphi(A)}$ is a fuzzy subalgebra of $L'$. Hence, for $y_1, y_2\in L'$ and $\alpha\in F$, we have
\begin{itemize}
\item[(i)] $r_{\varphi(A)}(y_1+y_2)\geq r_{\varphi(A)}(y_1)\wedge r_{\varphi(A)}(y_2)$,
\item[(ii)] $r_{\varphi(A)}(\alpha y_1))\geq r_{\varphi(A)}(y_1)$,
\item[(iii)] $r_{\varphi(A)}([y_1, y_2])\geq r_{\varphi(A)}(y_1)\wedge r_{\varphi(A)}(y_2)$.
\end{itemize}
Now our result follows from the homogeneity of $\varphi(A)$. \hfill $\Box$

Chung-Gook Kim and Dong-Soo Lee (\cite{Kim}) proved if $\varphi:L\rightarrow L'$ is a surjective Lie algebra homomorphism and $A=\mu_A$ is a fuzzy ideal of $L$, then $\varphi(A)$ is a fuzzy ideal of $L'$. This fact and Theorem \ref{CLF12A} help us to extend it to complex fuzzy Lie algebra case.
\begin{theorem}\label{CFL11}
Let $\varphi:L\rightarrow L'$ be a surjective Lie algebra homomorphism. If $A=\mu_A$, where $\mu_A(x)=r_A(x)e^{iw_A(x)}$ for any $x\in L$, is a complex fuzzy ideal of $L$, then $\varphi(A)$ is also a complex fuzzy ideal of $L'$.
\end{theorem}

\begin{theorem}\label{CFL12}
Let $\varphi: L\rightarrow L'$ be a surjective Lie homomorphism. If $A=\mu_A$ and $B=\mu_B$ are complex fuzzy ideals of $L$ such that $A$ is homogeneous of $B$, then $\varphi(A+B)= \varphi(A)+ \varphi(B)$.
\end{theorem}
{\it Proof.~}
For $y\in L'$, we have
\begin{eqnarray*}
\mu_{\varphi(A)+\varphi(B)}(y)&=&\mathrm{sup}_{y=\varphi(x)}\{\mu_{A+B}(x)\}\\
&=& \mathrm{sup}_{y=\varphi(x)}\{\mathrm{sup}_{x=a+b}\{\mu_A(a)\wedge\mu_B(b)\}\}\\
&=&\mathrm{sup}_{y=\varphi(a)+ \varphi(b)}\{\mu_A(a)\wedge\mu_B(b)\}\\
&=&\mathrm{sup}_{y=m+n}\{\mathrm{sup}_{m=\varphi(a)}\{\mu_A(a)\}\wedge\mathrm{sup}_{n=\varphi(b)}\{\mu_B(b)\}\}\\
&=&\mathrm{sup}_{y=m+n}\{\mu_{\varphi(A)}(m)\wedge\mu_{\varphi(B)}(n)\}\\
&=&\mu_{\varphi(A)+\varphi(B)}(y).
\end{eqnarray*}
\hfill $\Box$

Let $A=\mu_A$ be a complex fuzzy subset of a set $X$ with membership function $\mu_A(x)=r_A(x)e^{iw_A(x)}$. For $\alpha\in [0, 1]$ and $\beta \in [0, 2\pi]$, we define the sets:
\begin{itemize}
\item[(i)] $A_{\alpha, \beta}=\{x\in X~|~r_A(x)\geq \alpha ~\mathrm{and}~ w_A(x)\geq \beta\}$,
\item[(ii)] $A_{\alpha^>, \beta}=\{x\in X~|~r_A(x)> \alpha ~\mathrm{and}~ w_A(x)\geq \beta\}$,
\item[(iii)] $A_{\alpha, \beta^>}=\{x\in X~|~r_A(x)\geq \alpha ~\mathrm{and}~ w_A(x)> \beta\}$,
\item[(iv)] and $A_{\alpha^>, \beta^>}=\{x\in X~|~r_A(x)> \alpha ~\mathrm{and}~ w_A(x)> \beta\}$.
\end{itemize}
\begin{theorem}\label{CFL13}
Let $\varphi:L\rightarrow L'$ be a homomorphism Lie algebra and $B$ be a complex fuzzy Lie subalgebra (ideal) of $L'$ with membership function $\mu_B(x)=r_B(x)e^{iw_B(x)}$. Then, for $\alpha\in [0, 1]$ and $\beta \in [0, 2\pi]$, we have $\varphi^{-1}\left(B_{\alpha, \beta}\right)= \left(\varphi^{-1}(B)\right)_{\alpha, \beta}$.
\end{theorem}
{\it Proof.~}
$x\in \varphi^{-1}\left(B_{\alpha, \beta}\right)$ if and only if $\varphi(x)\in B_{\alpha, \beta}$ if and only if $\mu_B(\varphi(x))\geq \alpha e^{i\beta}$ if and only if $\mu_{\varphi^{-1}(B)}(x)\geq \alpha e^{i\beta}$ if and only if $x\in \left(\varphi^{-1}(B)\right)_{\alpha, \beta}$.\hfill $\Box$

In a similar way, we can prove the following theorems.
\begin{theorem}\label{CFL14}
Let $\varphi:L\rightarrow L'$ be a homomorphism Lie algebra and $B$ be a complex fuzzy Lie subalgebra (ideal) of $L'$ with membership function $\mu_B(x)=r_B(x)e^{iw_B(x)}$. Then, for $\alpha\in [0, 1]$ and $\beta \in [0, 2\pi]$, we have $\varphi^{-1}\left(B_{\alpha^>, \beta}\right)= \left(\varphi^{-1}(B)\right)_{\alpha^>, \beta}$.
\end{theorem}
\begin{theorem}\label{CFL15}
Let $\varphi:L\rightarrow L'$ be a homomorphism Lie algebra and $B$ be a complex fuzzy Lie subalgebra (ideal) of $L'$ with membership function $\mu_B(x)=r_B(x)e^{iw_B(x)}$. Then, for $\alpha\in [0, 1]$ and $\beta \in [0, 2\pi]$, we have $\varphi^{-1}\left(B_{\alpha, \beta^>}\right)= \left(\varphi^{-1}(B)\right)_{\alpha, \beta^>}$.
\end{theorem}
\begin{theorem}\label{CFL16}
Let $\varphi:L\rightarrow L'$ be a homomorphism Lie algebra and $B$ be a complex fuzzy Lie subalgebra (ideal) of $L'$ with membership function $\mu_B(x)=r_B(x)e^{iw_B(x)}$. Then, for $\alpha\in [0, 1]$ and $\beta \in [0, 2\pi]$, we have $\varphi^{-1}\left(B_{\alpha^>, \beta^>}\right)= \left(\varphi^{-1}(B)\right)_{\alpha^>, \beta^>}$.
\end{theorem}

\end{document}